\documentclass{svjour3}                     
\smartqed  
\usepackage{amstext,amssymb,amsmath}
\usepackage{graphicx,array,booktabs}
\newcommand{\sgn}{\mathop{\rm sgn}\nolimits}

\newcommand{\bfb}{\boldsymbol{b}}
\newcommand{\bfc}{\boldsymbol{c}}
\newcommand{\bfd}{\boldsymbol{d}}

\newcommand{\bfr}{\boldsymbol{r}}

\newcommand{\bfv}{\boldsymbol{v}}

\newcommand{\bfx}{\boldsymbol{x}}
\newcommand{\bfy}{\boldsymbol{y}}

\newcommand{\bfA}{\boldsymbol{A}}
\newcommand{\bfE}{\boldsymbol{E}}
\newcommand{\bfH}{\boldsymbol{H}}
\newcommand{\bfQ}{\boldsymbol{Q}}
\newcommand{\bfalpha}{\boldsymbol{\alpha}}
\newcommand{\bfbeta}{\boldsymbol{\beta}}
\newcommand{\bfgamma}{\boldsymbol{\gamma}}
\journalname{Under review}
\begin{document}

\title{MM Algorithms for Geometric and Signomial Programming\thanks{Research was supported by United States Public Health Service grants GM53275 and MH59490.}
}


\author{Kenneth Lange         \and
        Hua Zhou 
}


\institute{K. Lange \at
Departments of Biomathematics, Human Genetics, and Statistics, University of California, Los Angeles, CA 90095-1766, USA.	\\
\email{klange@ucla.edu}           
\and
H. Zhou \at
Department of Statistics, North Carolina State University, 2311 Stinson Drive, Campus Box 8203, Raleigh, NC 27695-8203, USA.	\\
\email{huazhou@ucla.edu}
}

\date{Received: date / Accepted: date}

\maketitle

\begin{abstract}
This paper derives new algorithms for signomial programming, a generalization of geometric 
programming. The algorithms are based on a generic principle for optimization called the MM 
algorithm. In this setting, one can apply the geometric-arithmetic mean inequality and a 
supporting hyperplane inequality to create a surrogate function with parameters separated.
Thus, unconstrained signomial programming reduces to a sequence of one-dimensional minimization 
problems.  Simple examples demonstrate that the MM algorithm derived can converge to a boundary point
or to one point of a continuum of minimum points.  Conditions under which the minimum point is 
unique or occurs in the interior of parameter space are proved for geometric programming. Convergence
to an interior point occurs at a linear rate. Finally, the MM framework easily accommodates
equality and inequality constraints of signomial type.  For the most important special case,
constrained quadratic programming, the MM algorithm involves very simple updates.
\end{abstract}

\section{Introduction}

As a branch of convex optimization theory, geometric programming is
next in line to linear and quadratic programming in importance
\cite{boyd04,ecker80,peressini88,peterson76}.  It has applications in chemical 
equilibrium problems \cite{passy68}, structural mechanics \cite{ecker80},
integrated circuit design \cite{hershenson01}, maximum likelihood estimation 
\cite{mazumdar83}, stochastic processes \cite{feigin81}, and a host of other 
subjects \cite{ecker80}. Geometric programming deals with posynomials, 
which are functions of the form 
\begin{eqnarray}
f(\bfx) & = & \sum_{\bfalpha \in S} c_{\bfalpha} \prod_{i=1}^n x_i^{\alpha_{i}} .
\label{general_posynomial}
\end{eqnarray}
Here the index set $S \subset \mathbb{R}^n$ is finite, and all 
coefficients $c_{\bfalpha}$ and all components $x_1,\ldots,x_n$ of the 
argument $\bfx$ of $f(\bfx)$ are positive.  The possibly fractional
powers $\alpha_i$ corresponding to a particular $\bfalpha$ may be positive, 
negative, or zero.  For instance, $x_1^{-1}+2x_1^3x_2^{-2}$ is a posynomial 
on $\mathbb{R}^2$.  In geometric programming we minimize 
a posynomial $f(\bfx)$ subject to posynomial inequality constraints of 
the form $u_j(\bfx) \le 1$ for $1 \le j \le q$, where the $u_j(\bfx)$ are again posynomials.
In some versions of geometric programming, equality constraints of
posynomial type are permitted \cite{boyd07}.

A signomial function has the same form as the posynomial (\ref{general_posynomial}), but the coefficients $c_{\bfalpha}$ are allowed to be negative. A signomial program is a generalization of a geometric program, where the objective and constraint functions can be signomials. From a computational point of view, signomial programming problems are significantly harder to solve than geometric programming problems. After suitable change of variables, a geometric program can be transformed into a convex optimization problem and globally solved by standard methods. In contrast, signomials may have many local minima.  Wang et al.\ \cite{WangZhangShen02SP} recently derived a path algorithm for solving unconstrained signomial programs.

The theory and practice of geometric programming has been stable for
a generation, so it is hard to imagine saying anything novel about 
either. The attractions of geometric programming include its beautiful
duality theory and its connections with the arithmetic-geometric
mean inequality.  The present paper derives new algorithms for
both geometric and signomial programming based on a generic device for iterative 
optimization called the MM algorithm \cite{hunter04,lange00}.  The MM perspective
possesses several advantages. First it provides a unified framework for solving both geometric and 
signomial programs.  The algorithms derived here operate by separating 
parameters and reducing minimization of the objective function to a 
sequence of one-dimensional minimization problems. Separation of parameters
is apt to be an advantage in high-dimensional problems.  Another advantage is ease of implementation 
compared to competing methods of unconstrained geometric and signomial programming \cite{WangZhangShen02SP}.
Finally, straightforward generalizations of our MM algorithms extend beyond signomial programming.

We conclude this introduction by sketching a roadmap to the rest of the paper. Section \ref{sec:MM} reviews the MM algorithm. Section \ref{sec:MM-unconst} derives MM algorithm for unconstrained signomial program from two simple inequalities. The behavior of the MM algorithm is illustrated on a few numerical examples in Section \ref{sec:unconst-examples}. Section \ref{sec:MM-constr} extends the MM algorithm for unconstrained problems to the constrained cases using the penalty method. Section \ref{sec:constr-examples} specializes to linearly constrained quadratic programming on the positive orthant. Convergence results are discussed in Section \ref{sec:convergence}. 

\section{Background on the MM Algorithm}
\label{sec:MM}

The MM principle involves majorizing the objective function $f(\bfx)$ 
by a surrogate function $g(\bfx \mid \bfx_m)$ around the current iterate $\bfx_m$ 
(with $i$th component $x_{mi}$) of a search.  Majorization is defined by the 
two conditions
\begin{eqnarray}
f(\bfx_m) & = & g(\bfx_m \mid \bfx_m)  \label{majorization_definition} \\
f(\bfx) & \le & g(\bfx \mid \bfx_m)\: , \quad \quad \bfx \ne \bfx_m . \nonumber
\end{eqnarray}
In other words, the surface $\bfx \mapsto g(\bfx \mid \bfx_m)$ lies above the
surface $\bfx \mapsto f(\bfx)$ and is tangent to it at the point $\bfx=\bfx_m$.  Construction
of the majorizing function $g(\bfx \mid \bfx_m)$ constitutes the first
M of the MM algorithm.

The second M of the algorithm minimizes the surrogate $g(\bfx \mid \bfx_m)$
rather than $f(\bfx)$.  If $\bfx_{m+1}$ denotes the minimizer of $g(\bfx \mid \bfx_m)$, 
then this action forces the descent property $f(\bfx_{m+1}) \le f(\bfx_m)$.  This 
fact follows from the inequalities
\begin{eqnarray*}
f(\bfx_{m+1}) \le g(\bfx_{m+1} \mid \bfx_m) \le g(\bfx_m \mid \bfx_m) = f(\bfx_m), 
\end{eqnarray*}
reflecting the definition of $\bfx_{m+1}$ and the tangency conditions (\ref{majorization_definition}).  The descent property lends 
the MM algorithm remarkable numerical stability.  Strictly speaking, it depends 
only on decreasing $g(\bfx \mid \bfx_m)$, not on minimizing $g(\bfx \mid \bfx_m)$. 

\section{Unconstrained Signomial Programming}
\label{sec:MM-unconst}

The art in devising an MM algorithm revolves around intelligent choice
of the majorizing function.  For signomial programming problems, fortunately one can invoke two simple inequalities. For terms with positive coefficients $c_{\bfalpha}$, we use the arithmetic-geometric mean inequality  
\begin{eqnarray}
\prod_{i=1}^n z_i^{\alpha_i} & \le & \sum_{i=1}^n \frac{\alpha_i}{\|\bfalpha\|_1} z_i^{\|\bfalpha\|_1}
\label{original-arithmetic-geometric}
\end{eqnarray}
for nonnegative numbers $z_i$ and $\alpha_i$ and $\ell_1$ norm $\|\bfalpha\|_1 = \sum_{i=1}^n |\alpha_i|$
\cite{steele04}.  If we make the choice  $z_i = x_i/x_{mi}$ in inequality
(\ref{original-arithmetic-geometric}), then the majorization
\begin{eqnarray}
\prod_{i=1}^n x_i^{\alpha_i} & \le & \left(\prod_{i=1}^n x_{mi}^{\alpha_i}\right)
\sum_{i=1}^n \frac{\alpha_i}{\|\bfalpha\|_1} \left( \frac{x_i}{x_{mi}}\right)^{\|\bfalpha\|_1}, 
\label{majorizing-arithmetic-geometric}
\end{eqnarray}
emerges, with equality when $\bfx= \bfx_m$.  We can broaden the scope of the 
majorization (\ref{majorizing-arithmetic-geometric}) to cases with $\alpha_i<0$ by
replacing $z_i$ by the reciprocal ratio $x_{mi}/x_i$ whenever $\alpha_i<0$. Thus, for terms $c_{\bfalpha} \prod_{i=1}^n x_i^{\alpha_{i}}$ with $c_{\bfalpha}>0$, we have the majorization
\begin{eqnarray*}
c_{\bfalpha} \prod_{i=1}^n x_i^{\alpha_i} \le c_{\bfalpha} \left( \prod_{j=1}^n x_{mj}^{\alpha_j} \right) \sum_{i=1}^n \frac{|\alpha_i|}{\|\bfalpha\|_1} \left( \frac{x_i}{x_{mi}} \right)^{\|\bfalpha\|_1 \text{sgn}(\alpha_i)},
\end{eqnarray*}
where $\sgn(\alpha_i)$ is the sign function.

The terms $c_{\bfalpha} \prod_{i=1}^n x_i^{\alpha_{i}}$
with $c_{\bfalpha} <0$ are handled by a different majorization. Our point of departure is
the supporting hyperplane minorization
\begin{eqnarray*}
z & \ge & 1 + \ln z 
\end{eqnarray*}
at the point $z=1$. If we let $z = \prod_{i=1}^n (x_i/x_{mi})^{\alpha_{i}}$,
then it follows that
\begin{eqnarray}
\prod_{i=1}^n x_i^{\alpha_{i}} & \ge & \prod_{j=1}^n x_{mj}^{\alpha_{j}}
\left(1 + \sum_{i=1}^n \alpha_i \ln x_i - \sum_{i=1}^n \alpha_i \ln x_{mi}\right) \label{eqn:supp-hp-maj}
\end{eqnarray}
is a valid minorization in $\bfx$ around the point $\bfx_m$.  Multiplication
by the negative coefficient $c_{\bfalpha}$ now gives the desired majorization. The surrogate
function separates parameters and is convex when all of the $\alpha_i$ are positive.

In summary, the objective function (\ref{general_posynomial}) is majorized up to an
irrelevant additive constant by the sum 
\begin{eqnarray}
g(\bfx \mid \bfx_m) & =  & \sum_{i=1}^n g_i(x_i \mid \bfx_m) \nonumber \\
g_i(x_i \mid \bfx_m) &=  & \sum_{\bfalpha \in S_+} c_{\bfalpha} \Bigg( \prod_{j=1}^n x_{mj}^{\alpha_j} \Bigg) \frac{|\alpha_i|}{\|\bfalpha\|_1} \left( \frac{x_i}{x_{mi}} \right)^{\|\bfalpha\|_1 \text{sgn}(\alpha_i)}	 \label{signomial_majorizer} \\
 &   &  + \sum_{\bfalpha \in S_-} c_{\bfalpha} \Bigg( \prod_{j=1}^n x_{mj}^{\alpha_j} \Bigg) \alpha_i \ln x_i, \nonumber
\end{eqnarray}
where $S_+=\{\bfalpha:c_{\bfalpha}>0\}$, and $S_-=\{\bfalpha:c_{\bfalpha}<0\}$.  To guarantee that the next iterate is well defined and occurs on the interior of the parameter domain, it is helpful to assume for each $i$ that at least one $\bfalpha \in S_{+}$ has $\alpha_i$ positive and at least one $\bfalpha \in S_{+}$ has $\alpha_i$ negative.  Under these conditions each $g_i(x_i \mid \bfx_m)$ is coercive and attains its minimum on the open interval $(0,\infty)$. 

Minimization of the majorizing function is straightforward because the surrogate functions $g_i(x_i \mid \bfx_m)$ are univariate functions. The derivative of $g_i(x_i \mid \bfx_m)$ with respect to its left argument equals
\begin{eqnarray*}
g_i'(x_i \mid \bfx_m) & = & \sum_{\bfalpha \in S_+} c_{\bfalpha} 
\Bigg(\prod_{j=1}^n x_{mj}^{\alpha_j}\Bigg)
\alpha_i x_i^{-1}
\left(\frac{x_i}{x_{mi}}\right)^{ \|\bfalpha\|_1 \sgn(\alpha_i)}	\\
& & \hspace{.5in} + \sum_{\bfalpha \in S_-} c_{\bfalpha} \Bigg( \prod_{j=1}^n x_{mj}^{\alpha_j} \Bigg) \alpha_i x_i^{-1}
\end{eqnarray*}
Assuming that the exponents $\alpha_i$ are integers, this is a rational function of $x_i$, and once we equate it to 0, we are faced with solving a polynomial equation.  This task can be accomplished by bisection or by Newton's method.  

In a geometric program, the function  $g_i'(x_i \mid \bfx_m)$ has a single root on the interval $(0,\infty)$.  For a proof of this fact,
note that making the standard change of variables $x_i = e^{y_i}$ eliminates the positivity constraint $x_i>0$ and renders the transformed function $h_i(y_i \mid \bfx_m)=g_i(x_i \mid \bfx_m)$ strictly convex.  Because $|\alpha_i| \sgn(\alpha_i)^2 = |\alpha_i|$, the 
second derivative 
\begin{eqnarray*}
h_i''(y_i \mid \bfx_m) & = & \sum_{\bfalpha \in S_+} c_{\bfalpha} 
\Bigg(\prod_{j=1}^n x_{mj}^{\alpha_j}\Bigg)
\frac{|\alpha_i| \cdot \|\bfalpha\|_1}{x_{mi}^{ \|\bfalpha\|_1 \sgn(\alpha_i)}} 
e^{\|\bfalpha\|_1 \sgn(\alpha_i)y_i}
\end{eqnarray*}
is positive.  Hence, $h_i(y_i \mid \bfx_m)$ is strictly convex and possesses 
a unique minimum point.  These arguments yield the even sweeter dividend that
the MM iteration map is continuously differentiable. From the
vantage point of the implicit function theorem \cite{hoffman75}, the stationary condition
$h_i'(y_{m+1,i} \mid \bfx_m)=0$ determines $y_{m+1,i}$, and consequently
$x_{m+1,i}$, in terms of $\bfx_m$. Observe here that $h_i''(y_{mi} \mid \bfx_m) \neq 0$
as required by the implicit function.

It is also worth pointing out that even more functions can be brought
under the umbrella of signomial programming.  For instance, 
majorization of the functions $- \ln f(\bfx)$ and $\ln f(\bfx)$ is possible for
any posynomial $f(\bfx)=\sum_{\bfalpha} c_{\bfalpha} \prod_{i=1}^n x_i^{\alpha_i}$.  In the first case, 
\begin{eqnarray}
- \ln f(\bfx) & \le & -\sum_{\bfalpha} \frac{a_{m \bfalpha}}{b_m}\Big[\sum_{i=1}^n \alpha_i \ln x_i
+ \ln \Big(\frac{c_{\bfalpha} b_m }{a_{m\bfalpha}}\Big)\Big] \label{log_complication1}
\end{eqnarray} 
holds for $a_{m \bfalpha} = c_{\bfalpha} \prod_{i=1}^n x_{mi}^{\alpha_i}$ and
$b_m  = \sum_{\bfalpha} a_{m\bfalpha}$ because Jensen's inequality applies
to the convex function $-\ln t$.  In the second case, the supporting hyperplane 
inequality applied to the convex function $-\ln t$ implies
\begin{eqnarray*}
\ln f(\bfx) & \le & \ln f(\bfx_m)+\frac{1}{f(\bfx_m)}\Big[f(\bfx)-f(\bfx_m)\Big].
\end{eqnarray*} 
This puts us back in the position of needing to majorize a posynomial, a
problem we have already discussed in detail. By our previous remarks, the 
coefficients  $c_{\bfalpha}$ can be negative as well as positive in this case.
Similar majorizations apply to any composition $\phi \circ f(\bfx)$ of
a posynomial $f(\bfx)$ with an arbitrary concave function $\phi(y)$.

\section{Examples of Unconstrained Minimization}
\label{sec:unconst-examples}

Our first examples demonstrate the robustness of the MM algorithms in minimization
and illustrate some of the complications that occur.  In each case we can explicitly calculate
the MM updates. To start, consider the posynomial
\begin{eqnarray*}
f_1(\bfx) & = & \frac{1}{x_1^3}+\frac{3}{x_1x_2^2}+x_1x_2
\end{eqnarray*}
with the implied constraints $x_1>0$ and $x_2>0$.  The majorization 
(\ref{majorizing-arithmetic-geometric}) applied to the
third term of $f_1(\bfx)$ yields
\begin{eqnarray*}
x_1x_2 & \le & x_{m1}x_{m2}\left[\frac 12 \left( \frac{x_1}{x_{m1}}\right)^2
+ \frac 12 \left( \frac{x_2}{x_{m2}} \right)^2 \right] \\
& = & \frac{x_{m2}}{2 x_{m1}} x_1^2+\frac{x_{m1}}{2 x_{m2}}x_2^2.
\end{eqnarray*}
Applied to the second term of $f_1(\bfx)$ using the reciprocal ratios, it gives 
\begin{eqnarray*}
{3 \over x_1x_2^2} & \le & {3 \over x_{m1}x_{m2}^2}
\left[{1 \over 3} \left( {x_{m1} \over x_1}\right)^3
+{2 \over 3}\left( {x_{m2} \over x_2} \right)^3 \right] \\
& = &  {x_{m1}^2 \over x_{m2}^2} {1 \over x_1^3}+{2 x_{m2} \over x_{m1}}{1 \over x_2^3}.
\end{eqnarray*}
The sum $g(\bfx \mid \bfx_m)$ of the two surrogate functions
\begin{eqnarray*}
g_1(x_1 \mid \bfx_m) & = & {1 \over x_1^3}+{x_{m1}^2 \over x_{m2}^2} {1 \over x_1^3}
+{x_{m2} \over 2 x_{m1}} x_1^2 \\
g_2(x_2 \mid \bfx_m) & = & {2 x_{m2} \over x_{m1}}{1 \over x_2^3}
+ {x_{m1} \over 2 x_{m2}}x_2^2
\end{eqnarray*}
majorizes $f_1(\bfx)$. If we set the derivatives
\begin{eqnarray*}
g_1'(x_1 \mid \bfx_m) & = & -{3 \over x_1^4}
-{x_{m1}^2 \over x_{m2}^2} {3 \over x_1^4}
+{x_{m2} \over x_{m1}} x_1 \\
g_2'(x_1 \mid \bfx_m) & = & -{6 x_{m2} \over x_{m1}}{1 \over x_2^4}
+ {x_{m1} \over x_{m2}}x_2
\end{eqnarray*}
of each of these equal to 0, then the updates
\begin{eqnarray*}
x_{m+1,1} & = &\sqrt[5]{3\left({x_{m1}^2 \over x_{m2}^2}+1\right){x_{m1} \over x_{m2}}}, \quad \quad
x_{m+1,2} \;\;\, = \;\;\, \sqrt[5]{6 {x_{m2}^2 \over x_{m1}^2}}
\end{eqnarray*}
solve the minimization step of the MM algorithm.  It is also obvious that the point 
$\bfx=(\sqrt[5]{6},\sqrt[5]{6})^t$ is a fixed point of the updates, and the reader
can check that it minimizes $f_1(\bfx)$.  


It is instructive to consider the slight variations
\begin{eqnarray*}
f_2(\bfx) & = & {1 \over x_1x_2^2}+x_1x_2^2 \\
f_3(\bfx) & = & {1 \over x_1x_2^2}+x_1x_2 
\end{eqnarray*}
of this objective function.  In the first case, the reader can check that the MM 
algorithm iterates according to
\begin{eqnarray*}
x_{m+1,1} & = & \sqrt[3]{{x_{m1}^2 \over x_{m2}^2}}, \quad \quad
x_{m+1,2} \;\;\, = \;\;\, \sqrt[3]{{x_{m2} \over x_{m1}}} .
\end{eqnarray*}
In the second case, it iterates according to
\begin{eqnarray*}
x_{m+1,1} & = & \sqrt[5]{{x_{m1}^3 \over x_{m2}^3}}, \quad \quad
x_{m+1,2} \;\;\, = \;\;\, \sqrt[5]{2 {x_{m2}^2 \over x_{m1}^2}} .
\end{eqnarray*}
The objective function $f_2(\bfx)$ attains its minimum value whenever $x_1x_2^2 = 1$.  
The MM algorithm for $f_2(\bfx)$ converges after a single iteration to the
value 2, but the converged point depends on the initial point 
$\bfx_{0}$.  The infimum of $f_3(\bfx)$ is 0.  This value 
is attained asymptotically by the MM algorithm, which satisfies the 
identities $x_{m1}x_{m2}^{3/2} = 2^{3/10}$ and $x_{m+1,2}=2^{2/25}x_{m2}$ 
for all $m \ge 1$.  These results imply that $x_{m1}$ tends to 0 and 
$x_{m2}$ to $\infty$ in such a manner that $f_3(\bfx_m)$ tends to 0. 
One could not hope for much better behavior of the MM algorithm in 
these two examples.

The function
\begin{eqnarray*}
f_4(\bfx) & = &  x_1^2x_2^2 - 2x_1x_2x_3x_4 + x_3^2x_4^2	\;\; = \;\;
(x_1x_2-x_3x_4)^2
\end{eqnarray*}
is a signomial but not a posynomial. The surrogate function (\ref{signomial_majorizer})
reduces to
\begin{eqnarray*}
g(\bfx \mid \bfx_m) & = & \frac{x_{m2}^2}{2 x_{m1}^2}x_1^4+
\frac{x_{m1}^2}{2 x_{m2}^2}x_2^4+\frac{x_{m4}^2}{2 x_{m3}^2}x_3^4+\frac{x_{m3}^2}{2 x_{m4}^2}x_4^4 \\
&  & - 2x_{m1}x_{m2}x_{m3}x_{m4}(\ln x_1+\ln_2+\ln x_3+\ln x_4)
\end{eqnarray*}
with all variables separated. The MM updates
\begin{eqnarray*}
x_{m+1,1} & = & \sqrt[4]{\frac{x_{m1}^3x_{m3}x_{m4}}{x_{m2}}}, \quad \quad
x_{m+1,2} \;\; = \;\; \sqrt[4]{\frac{x_{m2}^3x_{m3}x_{m4}}{x_{m1}}}  \\
x_{m+1,3} & = & \sqrt[4]{\frac{x_{m3}^3x_{m1}x_{m2}}{x_{m4}}}, \quad \quad
x_{m+1,4} \;\; = \;\; \sqrt[4]{\frac{x_{m4}^3x_{m1}x_{m2}}{x_{m3}}} 
\end{eqnarray*}
converge in a single iteration to a solution of $f_4(\bfx) =0$. Again
the limit depends on the initial point.

The function
\begin{eqnarray*}
f_5(\bfx) & = & x_1x_2+x_1x_3+x_2x_3-\ln(x_1+x_2+x_3)
\end{eqnarray*}
is more complicated than a signomial.  It also is unbounded because the
point $\bfx$ with components $x_1=m$ and $x_2 = x_3= 1/m$ satisfies
 $f_5(\bfx) = 2+m^{-2}- \ln(m+2/m)$.  According to the majorization 
(\ref{log_complication1}), an appropriate surrogate is
\begin{eqnarray*}
g(\bfx \mid \bfx_m) & = & \Big(\frac{x_{m2}}{2x_{m1}}+\frac{x_{m3}}{2x_{m1}}\Big)x_1^2
+ \Big(\frac{x_{m1}}{2x_{m2}}+\frac{x_{m3}}{2x_{m2}}\Big)x_2^2+ \Big(\frac{x_{m1}}{2x_{m3}}
+\frac{x_{m2}}{2x_{m3}}\Big)x_3^2 \\
&  & - \frac{x_{m1}}{x_{m1}+x_{m2}+x_{m3}} \ln x_1 
 -\frac{x_{m2}}{x_{m1}+x_{m2}+x_{m3}} \ln x_2 \\
&   & - \frac{x_{m3}}{x_{m1}+x_{m2}+x_{m3}} \ln x_3 
\end{eqnarray*}
up to an irrelevant constant. The MM updates are
\begin{eqnarray*}
x_{m+1,i} & = & \sqrt{\frac{x_{mi}^2}{ (\sum_{j \ne i} x_{mj})(x_{m1}+x_{m2}+x_{m3})}} .
\end{eqnarray*}
If the components of the initial point coincide, then the iterates converge in a single iteration
to the saddle point with all components equal to $1/\sqrt{6}$.  Otherwise, it appears that 
$f_5(\bfx_m)$ tends to $-\infty$.

The following objective functions
\begin{eqnarray*}
f_6(\bfx) & = & x_1^2x_2^6 + x_1^2x_2^4 - 2x_1^2x_2^3 - x_1^2x_2^2 + 5.25 x_1x_2^3	\\
& & -2 x_1^2x_2 + 4.5 x_1 x_2^2 + 3x_1^2 + 3x_1 x_2 - 12.75 x_1	\\
f_7(\bfx) & = & \sum_{i=1}^{10} x_i^4 + 2 \sum_{i=1}^{9} x_i^2 \sum_{j=i+1}^{10} x_j^2 
+ (10^{-5}-0.5) \sum_{i=1}^{10} x_i^2 \\
&	& - (2 \times 10^{-5}) \sum_{i=7}^{10} x_i  + \frac{1}{16}	\\
f_8(\bfx) & = & x_1x_3^2x_6^{-1}x_7^{-1}+x_1^2x_3^{-1}x_5^{-2}x_6^{-1}x_7 	\\
&  & + x_1^3x_2^2x_5^{-2}x_6^2 + x_2^{-1}x_4^{-1}x_6^2 + x_3x_5^3x_6^{-3}	\\
f_9(\bfx) & = & x_1x_4^2 + x_2x_3 + x_1x_2x_3x_4^2 + x_1^{-1}x_4^{-2} 	
\end{eqnarray*}
from the reference \cite{WangZhangShen02SP} are intended for numerical illustration. 
Table \ref{table:unconstr-examples} lists initial conditions, minimum points, minimum values, and number of iterations until convergence
under the MM algorithm. Convergence is declared when the relative change in the objective function is less than a pre-specified value 
$\epsilon$, in other words, when
\begin{eqnarray*}
\frac{f(\bfx_m)-f(\bfx_{m+1})}{|f(\bfx_m)|+1} & \le & \epsilon.
\end{eqnarray*}
Optimization of the univariate surrogate functions easily succumbs to Newton's method. The MM algorithm takes fewer iterations to converge than the path algorithm for all of the test functions mentioned in \cite{WangZhangShen02SP} except $f_6(\bfx)$.  Furthermore, the MM algorithm avoids calculation of the gradient and Hessian and requires no matrix decompositions or selection of tuning constants.

\begin{figure}[htpb]
\begin{center}
$$
\begin{array}{cc}
\includegraphics[width=2.3in]{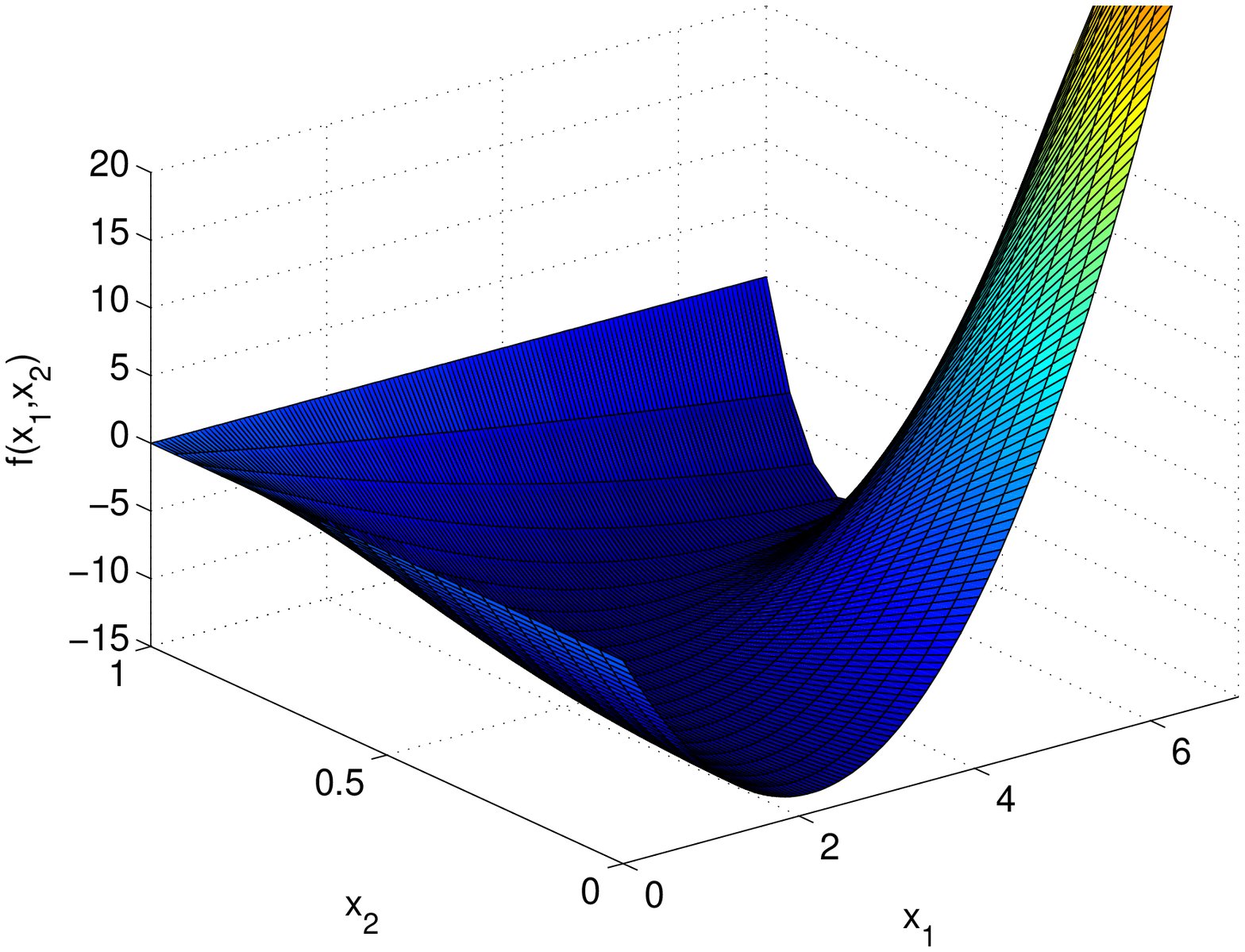} & \includegraphics[width=2.3in]{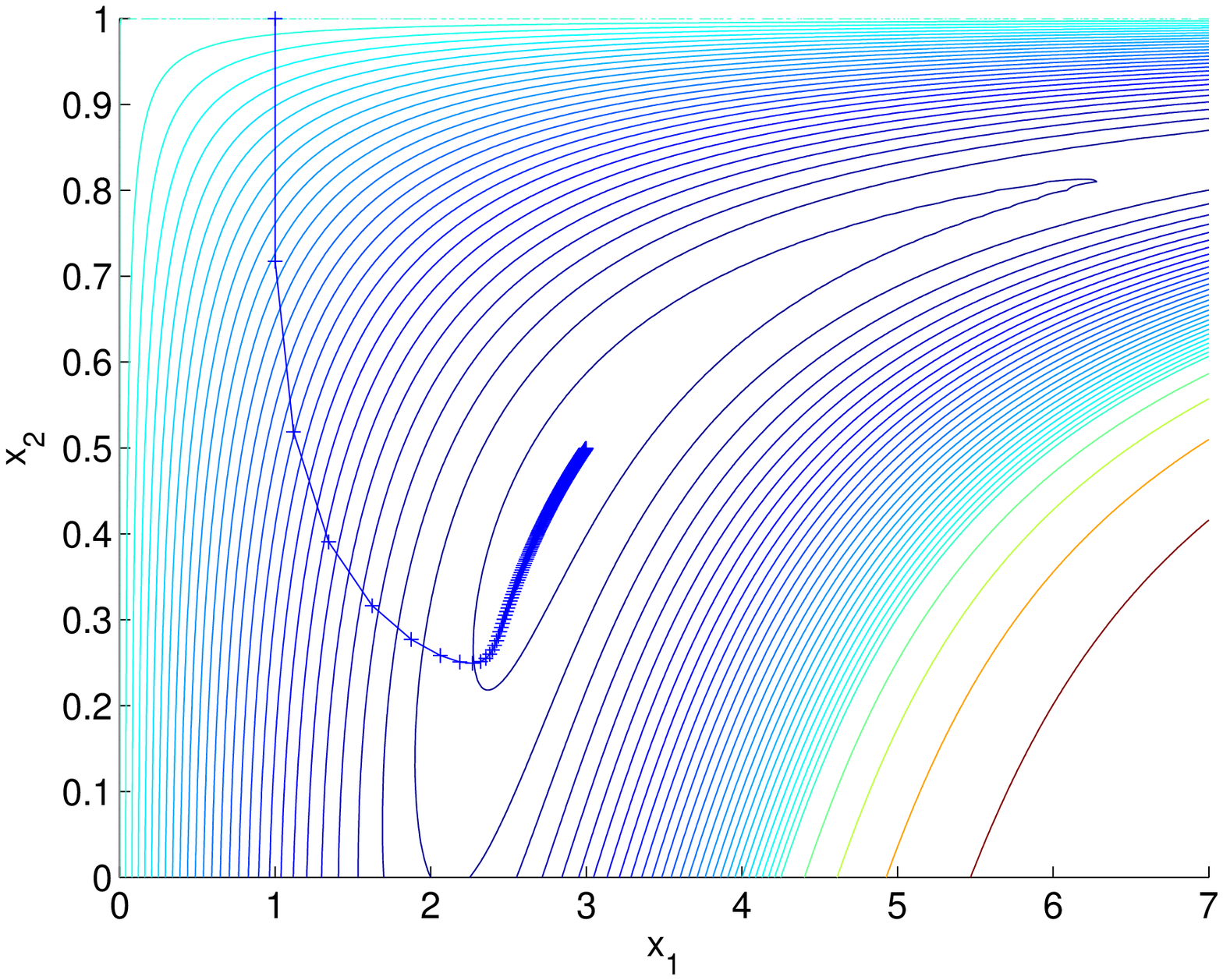}	\\
\includegraphics[width=2.3in]{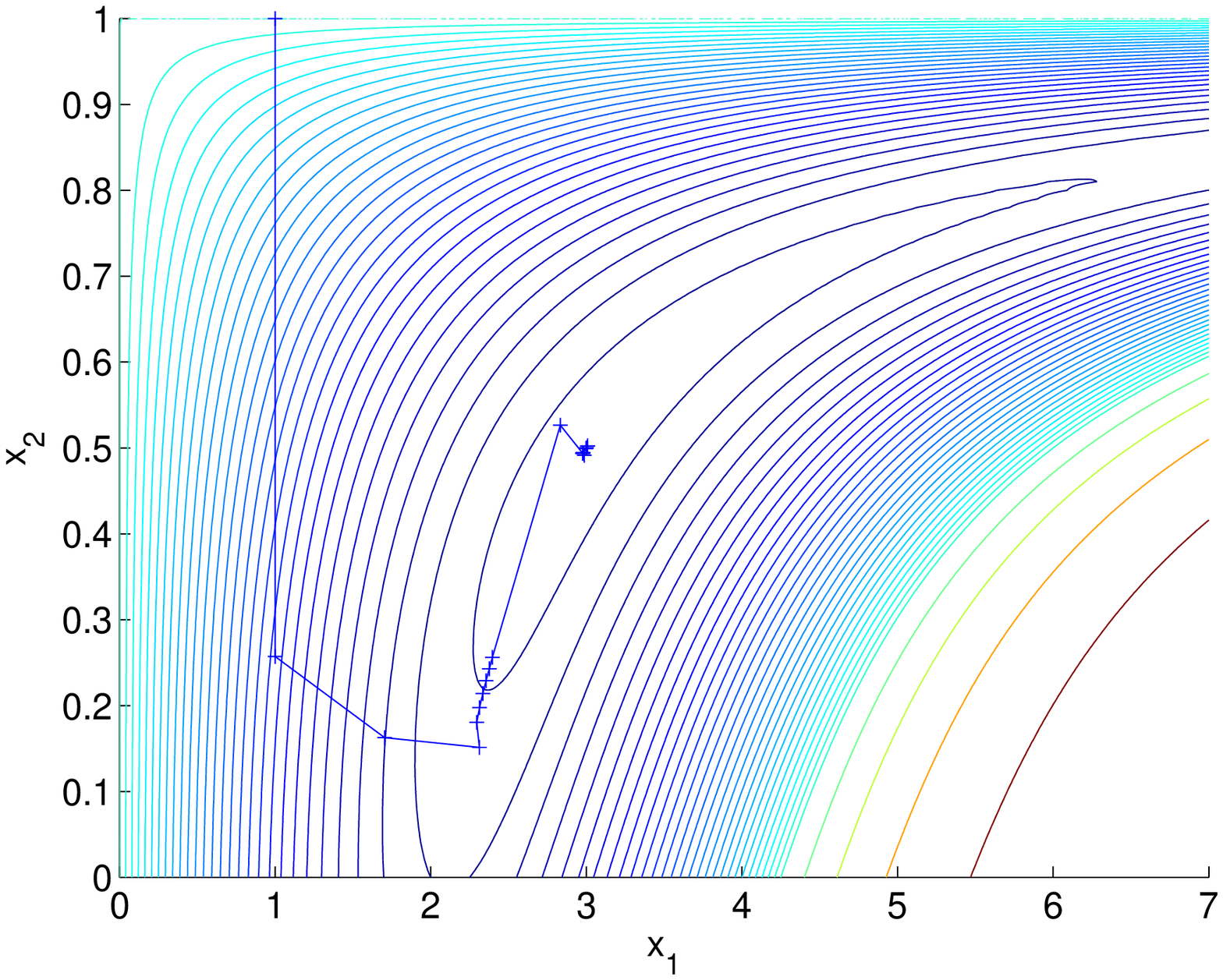} & \includegraphics[width=2.3in]{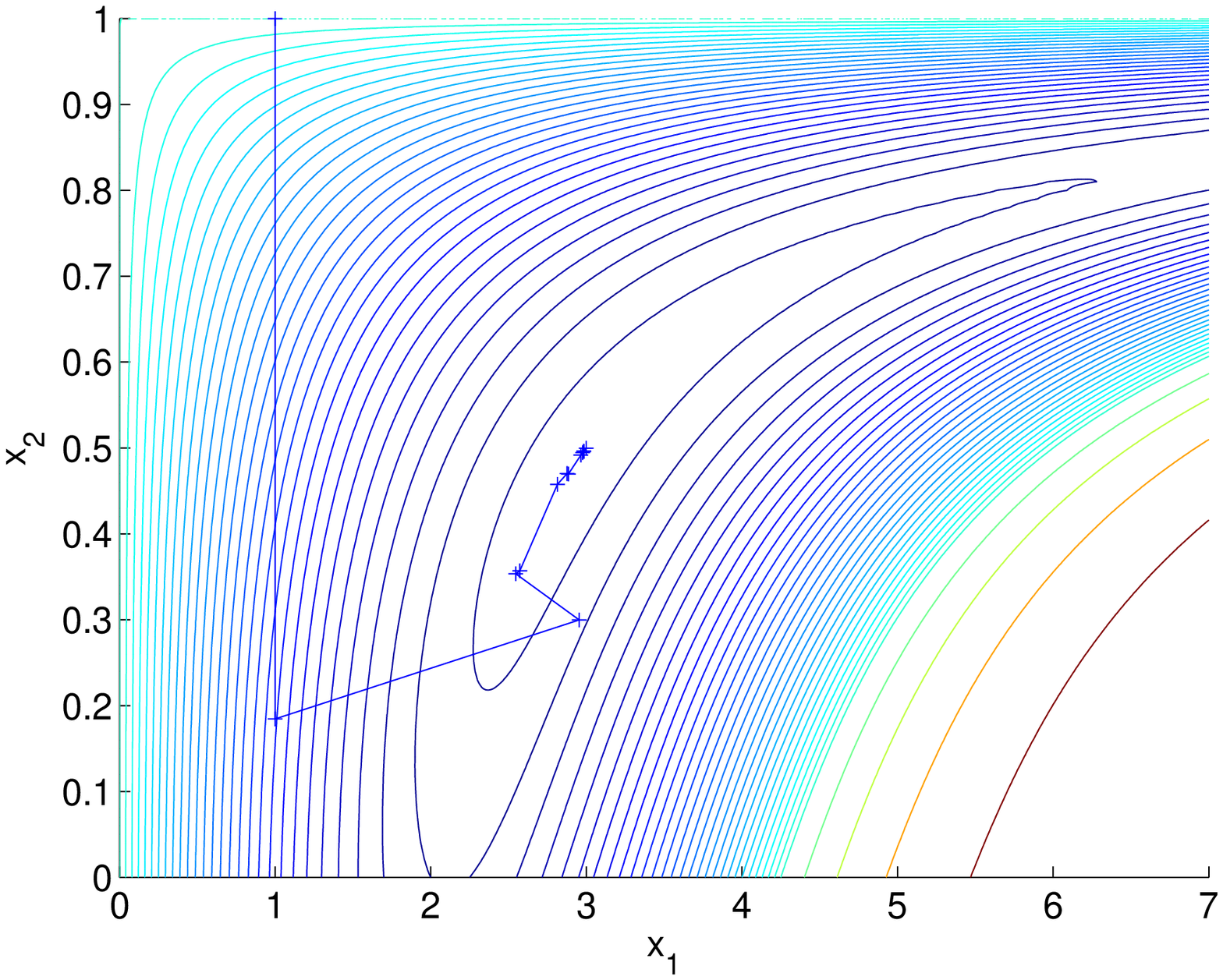}
\end{array}
$$
\end{center}
\caption{Upper left: The test function $f_6(\bfx)$. Upper right: 558 MM iterates. Lower left: 30 accelerated MM iterates ($q=1$ secant
conditions). Lower right: 12 accelerated MM iterates ($q=2$ secant conditions).}
\label{fig:test6-QN}
\end{figure}

As Section \ref{sec:convergence} observes, MM algorithms typically converge at a linear rate. Although slow convergence can occur for functions such as the test function $f_6(\bfx)$, there are several ways to accelerate an MM algorithm. For example, our published
quasi-Newton acceleration \cite{ZhouAlexanderLange09QN} often reduces the necessary number of iterations by one or two orders of magnitude. Figure \ref{fig:test6-QN} shows the progress of the MM iterates for the test function $f_6(\bfx)$ with and without quasi-Newton acceleration. Under a convergence criterion of $\epsilon=10^{-9}$ and $q=1$ secant condition, the required number of iterations falls to 30; under the same convergence criterion and $q=2$ secant conditions, the required number of iterations falls to 12. It is also worth
emphasizing that separation of parameters enables parallel processing in high-dimensional problems.  We have recently argued
\cite{ZhouLangeSuchard09MM-GPU} that the best approach to parallel processing is through graphics processing units (GPUs). These
cheap hardware devices offer one to two orders of magnitude acceleration in many MM algorithms with parameters separated. 

\begin{table}
\centering
{\scriptsize
\begin{tabular}{ccccccc}
\toprule
Fun & Type & Initial Point $\bfx_0$ & Min Point & Min Value & Iters ($10^{-9}$)	\\
\midrule
$f_1$ & P & (1,2) & (1.4310,1.4310) & 3.4128 & 38	\\
$f_2$ & P & (1,2) & (0.6300,1.2599) & 2.0000 & 2	\\
$f_3$ & P & (1,1) & diverges & 0.0000 &	\\
$f_4$ & S & (0.1,0.2,0.3,0.4) & (0.1596,0.3191,0.1954,0.2606) & 0.0000 & 3 	\\
$f_5$ & G & (1,1,1) & (0.4082,0.4082,0.4082) & 0.2973 & 2	\\
& & (1,2,3) & diverges & $-\infty$ & 	\\
$f_6$ & S & (1,1) & (2.9978,0.4994) & -14.2031 & 558 \\
$f_7$ & S & $(1,\ldots,10)$ & $0.0255 \bfx_0$ & 0.0000 & 18	\\
$f_8$ & P & $(1,\ldots,7)$ & diverges & 0.0000 &	\\
$f_9$ & P & (1,2,3,4) & (0.3969,0.0000,0.0000,1.5874) & 2.0000 & 7	\\
\bottomrule
\end{tabular}
}
\caption{Numerical examples of unconstrained signomial programming. Test functions $f_4(\bfx)$, $f_6(\bfx)$, $f_7(\bfx)$, $f_8(\bfx)$ and $f_9(\bfx)$ are taken from \cite{WangZhangShen02SP}. P: posynomial; S: signomial; G: general function.}
\label{table:unconstr-examples}
\end{table}

\section{Constrained Signomial Programming}
\label{sec:MM-constr}

Extending the MM algorithm to constrained geometric and signomial programming is challenging.  
Box constraints $a_i \le x_i \le b_i$ are consistent with parameter separation as just 
developed, but more complicated posynomial constraints that couple parameters are not.  
Posynomial inequality constraints take the form
\begin{eqnarray*}
h(\bfx) & = & \sum_{\bfbeta} d_{\bfbeta} \prod_{i=1}^n x_i^{\beta_{i}} 
\;\; \le \;\; 1.
\end{eqnarray*}
The corresponding equality constraint sets $h(\bfx)=1$.  We propose handling both constraints by penalty methods. Before we treat these matters in more depth, let us relax the positivity restrictions on the $d_{\bfbeta}$ but enforce the restriction $\beta_i \ge 0$.  The latter objective can be achieved by multiplying $h(\bfx)$ by $x_i^{\max_{\bfbeta} \{-\beta_i,0\}}$ for all $i$.  If we subtract the two sides of the resulting equality, then the equality constraint $h(\bfx)=1$ can be rephrased as $r(\bfx) = \sum_{\bfgamma} e_{\bfgamma} \prod_{i=1}^n x_i^{\gamma_{i}}=0$, with no restriction on the signs of the $e_{\bfgamma}$ but with the requirement $\gamma_i \ge 0$ in effect.  For example, the
equality constraint
\begin{eqnarray*}
\frac{1}{x_1}+\frac{x_1}{x_2^2} & = & 1
\end{eqnarray*}
becomes
\begin{eqnarray*}
x_1^2+x_2^2-x_1x_2^2 & = & 0.
\end{eqnarray*}

In the quadratic penalty method \cite{NocedalWrightBook,RuszczynskiBook} with objective function $f(\bfx)$ and a single equality constraint $r(\bfx)=0$ and a single inequality constraint $s(\bfx) \le 0$, one minimizes the sum $f_{\lambda}(\bfx) = f(\bfx)+\lambda r(\bfx)^2 + \lambda s(\bfx)_+^2$, where $s(\bfx)_+ = \max\{s(\bfx),0\}$.  As the penalty constant $\lambda$ tends to $\infty$, the solution vector $\bfx_{\lambda}$ typically converges to the constrained minimum. In the revised objective function, the term $r(\bfx)^2$ is a signomial whenever $r(\bfx)$ is a signomial. For example, in our toy problem  the choice $r(\bfx)= x_1^2+x_2^2-x_1x_2^2$ has square
\begin{eqnarray*}
r(\bfx)^2 & = & x_1^4+x_2^4+x_1^2x_2^4+2x_1^2x_2^2-2x_1x_2^4-2x_1^3x_2^2 .
\end{eqnarray*}
Of course, the powers in $r(\bfx)$ can be fractional here as well as integer. The term $s(\bfx)_+^2$ is not a signomial and must be subjected to the majorization
\begin{eqnarray*}
	s(\bfx)_+^2 \le \begin{cases}
	[s(\bfx)-s(\bfx_m)]^2	& s(\bfx_m) < 0	\\
	s(\bfx)^2 & s(\bfx_m) \ge 0
	\end{cases}
\end{eqnarray*}
to achieve this status.  In practice, one does not need to fully minimize $f_{\lambda}(\bfx)$  for any fixed $\lambda$. If one increases 
$\lambda$ slowly enough, then it usually suffices to merely decrease $f_{\lambda}(\bfx)$ at each iteration. The MM algorithm is designed to achieve precisely this goal. Our exposition so far suggests that we majorize $r(\bfx)^2$, $s(\bfx)^2$, and  $[s(\bfx)-s(\bfx_m)]^2$ in exactly the same manner that we majorize $f(\bfx)$. Separation of parameters generalizes, and the resulting MM algorithm keeps all parameters positive while permitting pertinent parameters to converge to 0.  Section \ref{sec:convergence} summarizes some of the convergence properties of this hybrid procedure.  

The quadratic penalty method traditionally relies on Newton's method to minimize the unconstrained functions $f_{\lambda}(\bfx)$.  Unfortunately, this tactic suffers from roundoff errors and numerical instability.  Some of these problems disappear with the MM algorithm.  No matrix inversions are involved, and iterates enjoy the descent property.  Ill-conditioning does cause harm in the form of slow convergence, but the previously mentioned quasi-Newton acceleration largely remedies the situation \cite{ZhouAlexanderLange09QN}.  As an alternative to quadratic penalties, exact penalties take the form $\lambda |r(\bfx)| + \lambda s(\bfx)_+$.  Remarkably, the exact penalty method produces the constrained minimum, not just in the limit, but for all finite $\lambda$ beyond a certain point.  Although this desirable property avoids the numerical instability encountered in the quadratic penalty method, the kinks in the objective functions  
$f(\bfx)+\lambda |r(\bfx)| + \lambda s(\bfx)_+$ are a nuisance. We will demonstrate in a future paper how to harness the MM algorithm to exact penalization.

\section{Nonnegative Quadratic Programming}
\label{sec:constr-examples}

As an illustration of constrained signomial programming,  consider quadratic programming over the positive orthant.  Let
\begin{eqnarray*}
f(\bfx) & = &  \frac{1}{2}  \bfx^t \bfQ \bfx + {\bf c}^t \bfx
\end{eqnarray*}
be the objective function, $\bfE \bfx = \bfd$ the linear equality constraints, and $\bfA \bfx \le \bfb$ the linear inequality constraints.   The symmetric matrix $\bfQ$ can be negative definite, indefinite, or positive definite. The quadratic penalty method involves minimizing the sequence of penalized objective functions
\begin{eqnarray*}
f_{\lambda}(\bfx) & = & \frac{1}{2} \bfx^t \bfQ \bfx + \bfc^t \bfx + \frac{\lambda}{2} \|(\bfA\bfx - \bfb)_+ \|_2^2 + \frac{\lambda}{2} \|\bfE\bfx - \bfd\|_2^2 
\end{eqnarray*}
as $\lambda$ tends to $\infty$. Based on the obvious majorization
\begin{eqnarray*}
x_+^2 & \le &
 \begin{cases}
(x-x_m)^2	& x_m < 0	\\
x^2 & x_m \ge 0
\end{cases},
\end{eqnarray*}
the term $\|(\bfA\bfx - \bfb)_+\|_2^2$ is majorized by $\|\bfA\bfx - \bfb - \bfr_{m}\|_2^2$, where 
\begin{eqnarray*}
\bfr_m & = & \min \{\bfA\bfx_m - \bfb, \bf0\}.
\end{eqnarray*}
A brief calculation shows that $f_{\lambda}(\bfx)$ is majorized by the surrogate function
\begin{eqnarray*}
g_{\lambda}(\bfx \mid \bfx_m) & = & \frac{1}{2} \bfx^t \bfH_{\lambda} \bfx + \bfv_{\lambda m}^t \bfx
\end{eqnarray*}
up to an irrelevant constant, where $\bfH_{\lambda}$ and $\bfv_{\lambda m}$ are defined by 
\begin{eqnarray*}
\bfH_\lambda  & = &  \bfQ + \lambda (\bfA^t \bfA + \bfE^t \bfE)	\\
\bfv_{\lambda m} & = & \bfc - \lambda \bfA^t (\bfb+\bfr_m) - \lambda \bfE^t \bfd.
\end{eqnarray*}
It is convenient to assume that the diagonal coefficients $\frac{1}{2}h_{\lambda ii}$ appearing in the quadratic form $\frac 12 \bfx^T \bfH_{\lambda} \bfx$ are positive.  This is generally the case for large $\lambda$.  One can handle the off-diagonal term $h_{\lambda ij} x_i x_j$  by either the majorization (\ref{majorizing-arithmetic-geometric}) or the majorization (\ref{eqn:supp-hp-maj}) according to the sign of $h_{\lambda ij}$. The reader can check that the MM updates reduce to 
\begin{eqnarray}
x_{m+1,i} & = & \frac{x_{mi}}{2} \left[- \frac{v_{\lambda mi}}{h_{\lambda mi}^+} + \sqrt{ \left(\frac{v_{\lambda mi}}{h_{\lambda mi}^+}\right)^2 - 4\frac{h_{\lambda mi}^-}{h_{\lambda mi}^+}} \, \right] ,	\label{eqn:qp-update2}
\end{eqnarray}
where
\begin{eqnarray*}
h_{\lambda mi}^+ &=& \sum_{j:h_{\lambda ij}>0} h_{\lambda ij} x_{mj}, \quad	 \quad h_{\lambda mi}^- = \sum_{j:h_{\lambda ij}<0} h_{\lambda ij} x_{mj}	. 
\end{eqnarray*}
When  $h_{\lambda mi}^- = 0$, the update (\ref{eqn:qp-update2}) collapses to
\begin{eqnarray}
x_{m+1,i} &=&  x_{mi} \max\Big\{- \frac{v_{\lambda mi}}{h_{\lambda mi}^+},0\Big\}. \label{eqn:qp-update3}
\end{eqnarray} 
To avoid sticky boundaries, we replace 0 in equation (\ref{eqn:qp-update3})  by a small positive constant $\epsilon$ such as $10^{-9}$. Sha et al.\ \cite{Sha} derived the update (\ref{eqn:qp-update2}) for $\lambda=0$ ignoring the constraints $\bfE \bfx = \bfd$ and $\bfA \bfx \le \bfb$.

For a numerical example without equality constraints take
\begin{eqnarray*}
f_{10}(\bfx) & = &  \frac{1}{2} x_1^2 + x_2^2 - x_1 x_2 - 2x_1 - 6x_2	\\
 \bfA & = &   \begin{pmatrix} 
 1 & 1 \\
 -1 & 2 \\
 2 & 1
 \end{pmatrix} , \quad \quad \bfb  \;\;\, =  \;\;\,\begin{pmatrix}
 2 \\
 2 \\
 3
 \end{pmatrix} .
\end{eqnarray*}
The minimum occurs at the point $(2/3,4/3)^t$. Table \ref{table:dummy2-quadpen} lists the number of iterations until convergence and the converged point $\bfx_\lambda$ for the sequence of penalty constants $\lambda=2^k$. The quadratic program
\begin{eqnarray*}
 f_{11}(\bfx) & = & -8x_1 - 16x_2 + x_1^2 + 4x_2^2	\\
  \bfA & = &   \begin{pmatrix} 
 1 & 1 \\
 1 & 0 
 \end{pmatrix} , \quad \quad \bfb  \;\;\, =  \;\;\,\begin{pmatrix}
 4 \\
 3
 \end{pmatrix} 
\end{eqnarray*}
converges much more slowly. Its minimum occurs at the point $(2.4,1.6)^t$. 
Table \ref{table:dummy-quadpen} lists the numbers of iterations until convergence with $(q=1$) and without ($q=0$) acceleration and the converged point $\bfx_\lambda$ for the same sequence of penalty constants $\lambda=2^k$.  Fortunately, quasi-Newton acceleration compensates for ill conditioning in this test problem.

\begin{table}
\centering
\begin{tabular}{r@{\hspace{.2in}}r@{\hspace{.2in}}c}
\toprule
$\log_2{\lambda}$ \hspace{-.1in} & Iters \hspace{-.1in} & $\bfx_\lambda$	\\
\midrule
0 & 8 & (0.9503,1.6464) \\   

1 & 6 & (0.8580,1.5164) \\   

2 & 5 & (0.8138,1.4461) \\   

3 & 23 & (0.7853,1.4067) \\   

4 & 32 & (0.7264,1.3702) \\   

5 & 31 & (0.6967,1.3518) \\   

6 & 30 & (0.6817,1.3426) \\   

7 & 29 & (0.6742,1.3380) \\   

8 & 28 & (0.6704,1.3356) \\   

9 & 26 & (0.6686,1.3345) \\   

10 & 25 & (0.6676,1.3339) \\   

11 & 23 & (0.6671,1.3336) \\   

12 & 22 & (0.6669,1.3335) \\   

13 & 21 & (0.6668,1.3334) \\   

14 & 19 & (0.6667,1.3334) \\   

15 & 18 & (0.6667,1.3334) \\   

16 & 16 & (0.6667,1.3333) \\   

17 & 15 & (0.6667,1.3333) \\   

%
%
%
\bottomrule
\end{tabular}
\caption{Iterates from the quadratic penalty method for the test function $f_{10}(\bfx)$. The convergence criterion for the inner loops is 
$10^{-9}$.}
\label{table:dummy2-quadpen}
\end{table}

\begin{table}
\centering
\begin{tabular}{r@{\hspace{.2in}}r@{\hspace{.3in}}r@{\hspace{.4in}}c}
\toprule
$\log_2{\lambda}$\hspace{-.1in} & Iters ($q=0$) \hspace{-.2in} & Iters ($q=1$)\hspace{-.3in} & $\bfx_\lambda$	\\
\midrule
0 & 18 & 5 & (3.0000,1.8000) \\   

1 & 2 & 2 & (2.8571,1.7143) \\   

2 & 56 & 6 & (2.6667,1.6667) \\   

3 & 97 & 5 & (2.5455,1.6364) \\   

4 & 167 & 5 & (2.4762,1.6190) \\   

5 & 312 & 5 & (2.4390,1.6098) \\   

6 & 541 & 6 & (2.4198,1.6049) \\   

7 & 955 & 5 & (2.4099,1.6025) \\   

8 & 1674 & 4 & (2.4050,1.6012) \\   

9 & 2924 & 3 & (2.4025,1.6006) \\   

10 & 4839 & 3 & (2.4013,1.6003) \\   

11 & 7959 & 4 & (2.4006,1.6002) \\   

12 & 12220 & 4 & (2.4003,1.6001) \\   

13 & 17674 & 4 & (2.4002,1.6000) \\   

14 & 21739 & 3 & (2.4001,1.6000) \\   

15 & 20736 & 3 & (2.4000,1.6000) \\   

16 & 8073 & 3 & (2.4000,1.6000) \\   

17 & 111 & 3 & (2.4000,1.6000) \\   

18 & 6 & 4 & (2.4000,1.6000) \\   

19 & 5 & 2 & (2.4000,1.6000) \\   

20 & 3 & 2 & (2.4000,1.6000) \\   

21 & 2 & 2 & (2.4000,1.6000) \\   
\bottomrule
\end{tabular}
\caption{Iterates from the quadratic penalty method for the test function $f_{11}(\bfx)$. The convergence criterion for the inner loops is 
$10^{-16}$.}
\label{table:dummy-quadpen}
\end{table}

\section{Convergence}
\label{sec:convergence}

As we have seen, the behavior of the MM algorithm is intimately
tied to the behavior of the objective function $f(\bfx)$.  For the sake
of simplicity, we now restrict attention to unconstrained minimization of
posynomials and investigate conditions guaranteeing that $f(\bfx)$ possesses
a unique minimum on its domain. Uniqueness is related to the 
strict convexity of the reparameterization
\begin{eqnarray*}
h(\bfy) & = & \sum_{\bfalpha \in S} c_{\bfalpha} e^{\bfalpha^t \bfy}
\end{eqnarray*}
of $f(\bfx)$, where $\bfalpha^t \bfy = \sum_{i=1}^n \alpha_i y_i$ is the
inner product of $\bfalpha$ and $\bfy$ and $x_i=e^{y_i}$ for each $i$.  The Hessian matrix
\begin{eqnarray*}
d^2 h(\bfy) & = & \sum_{\bfalpha \in S} c_{\bfalpha} e^{\bfalpha^t \bfy} \bfalpha \bfalpha^t
\end{eqnarray*}
of $h(\bfy)$ is positive semidefinite, so $h(\bfy)$ is convex.  If we let $T$ be the
subspace of $\mathbb{R}^n$ spanned by $\{\bfalpha\}_{\bfalpha \in S}$, then
$h(\bfy)$ is strictly convex if and only if $T = \mathbb{R}^n$.  Indeed, 
suppose the condition holds. For any $\bfv \neq {\bf 0}$, we then must have 
$\bfalpha^t \bfv \neq 0$ for some $\bfalpha \in S$.  It follows that 
\begin{eqnarray*}
\bfv^t d^2 h(\bfy)\bfv \;\; = \;\; \sum_{\bfalpha \in S} c_{\bfalpha} e^{\bfalpha^t \bfy} (\bfalpha^t \bfv)^2
\;\; > \;\; 0,
\end{eqnarray*}
and $d^2h(\bfy)$ is positive definite.  Conversely, suppose $T \neq \mathbb{R}^n$,
and take $\bfv \neq {\bf 0}$ with $\bfalpha^t \bfv = 0$ for every $\bfalpha \in S$.
Then $h(\bfy+t\bfv) = h(\bfy)$ for every scalar $t$, which is incompatible with 
$h(\bfy)$ being strictly convex.

Strict convexity guarantees uniqueness, not existence, of a minimum point.
Coerciveness ensures existence.  The objective function $f(\bfx)$ is coercive if
$f(\bfx)$ tends to $\infty$ whenever any component of $\bfx$ tends to 0 or
$\infty$.  Under the reparameterization $x_i=e^{y_i}$, this is equivalent
to $h(\bfy) = f(\bfx)$ tending to $\infty$ as $\|\bfy\|_2$ tends to $\infty$.
A necessary and sufficient condition for this to occur is that
$\max_{\bfalpha \in S} \bfalpha^t \bfv > 0$ for every $\bfv \neq {\bf 0}$.  For a
proof, suppose the contrary condition holds for some $\bfv \neq {\bf 0}$.  
Then it is clear that $h(t\bfv)$ remains bounded above by $h({\bf 0})$ as the 
scalar $t$ tends to $\infty$.  Conversely, if the stated condition 
is true, then the function $q(\bfy) = \max_{\bfalpha \in S} \bfalpha^t\bfy$ 
is continuous and achieves its minimum of $d>0$ on the sphere 
$\{\bfy \in \mathbb{R}^n : \|\bfy\|_2=1\}$.  It follows that $q(\bfy) \ge d\|\bfy\|_2$ and 
that
\begin{eqnarray*}
h(\bfy) \;\; \ge \;\; \max_{\bfalpha \in S} \{ c_{\bfalpha} e^{\bfalpha^ty}\}
\;\; \ge \;\; \left(\min_{\bfalpha \in S} c_{\bfalpha} \right) e^{d \|\bfy\|_2}.
\end{eqnarray*}
This lower bound shows that $h(\bfy)$ is coercive.

The coerciveness condition is hard to apply in practice.  An equivalent
condition is that the origin ${\bf 0}$ belongs to the interior of the
convex hull of the set $\{\bfalpha\}_{\bfalpha \in S}$.  It is straightforward
to show that the negations of these two conditions are logically
equivalent.  Thus, suppose $q(\bfv) = \max_{\alpha \in S} \bfalpha^t \bfv \le 0$ for some $\bfv \neq {\bf 0}$.
Every convex combination $\sum_{\bfalpha} p_{\bfalpha} \bfalpha$ then
satisfies $\left(\sum_{\bfalpha} p_{\bfalpha} \bfalpha \right)^t \bfv \le 0$.
If the origin is in the interior of the convex hull, then $\epsilon \bfv$ is 
also for every sufficiently small $\epsilon > 0$. But this leads 
to the contradiction $\epsilon \bfv^t \bfv = \epsilon \|\bfv\|_2^2 \le 0$.
Conversely, suppose ${\bf 0}$ is not in the interior of the
convex hull.  According to the separating hyperplane theorem for convex
sets, there exists a unit vector $\bfv$ with $\bfv^t \bfalpha \le 0 = \bfv^t{\bf 0}$
for every $\bfalpha \in S$.
In other words, $q(\bfv) \le 0$.  The convex hull criterion is easier to check, 
but it is not constructive.  In simple cases such as the objective
function $f_1(\bfx)$ where the power vectors are $\bfalpha=(-3,0)^t$, 
$\bfalpha=(-1,-2)^t$, and $\bfalpha=(1,1)^t$, it is visually obvious that
the origin is in the interior of their convex hull.

One can also check the criterion $q(\bfv) > 0$
for all $\bfv \neq {\bf 0}$ by solving a related geometric programming problem. 
This problem consists in minimizing the scalar $t$ subject to the inequality 
constraints $\bfalpha^t \bfy \le t$ for all $\bfalpha \in S$ and the nonlinear 
equality constraint $\|\bfy\|_2^2 = 1$.  If $t_{\mbox{\scriptsize min}} \le 0$, then 
the original criterion fails.

In some cases, the objective function $f(\bfx)$ does not attain its minimum on 
the open domain $\mathbb{R}_{>0}^n =\{\bfx: x_i>0, 1 \le i \le n\}$.  This condition is equivalent 
to the corresponding function $\ln h(\bfy)$ being unbounded below on 
$\mathbb{R}^n$.  According to Gordon's theorem \cite{borwein00,lange04}, 
this can happen if and only if ${\bf 0}$ is not in the convex hull of the
set $\{\bfalpha\}_{\bfalpha \in S}$.  Alternatively, both conditions
are equivalent to the existence of a vector $\bfv$ with $\bfalpha^t \bfv < 0$ 
for all $\bfalpha \in S$.  For the objective function $f_3(\bfx)$,
the power vectors are $\bfalpha=(-1,-2)^t$ and $\bfalpha=(1,1)^t$.  The 
origin $(0,0)^t$ does not lie on the line segment between them, and the
vector $(-3/2,1)^t$ forms a strictly oblique angle with each.  As
predicted, $f_3(\bfx)$ does not attain its infimum on $\mathbb{R}_{>0}^n$.

The theoretical development in reference \cite{lange04} 
demonstrates that the MM algorithm converges at a linear rate to the unique 
minimum point of the objective function $f(\bfx)$ when $f(\bfx)$ is coercive and 
its convex reparameterization $h(\bfy)$ is strictly convex. The theory does not
cover other cases, and it would be interesting to investigate them.
The general convergence theory of MM algorithms \cite{lange04} states that five properties of the objective function $f(\bfx)$ and MM algorithmic map $\bfx \mapsto M(\bfx)$ guarantee convergence to a stationary point of $f(\bfx)$: (a) $f(\bfx)$ is coercive on its open domain; (b) $f(\bfx)$ has only isolated stationary points; (c) $M(\bfx)$ is continuous; (d) $\bfx^*$ is a fixed point of $M(\bfx)$ if and only if $\bfx^*$ is a stationary point of $f(\bfx)$; and (e) $f[M(\bfx^*)] \ge f(\bfx^*)$, with equality if and only if $\bfx^*$ is a fixed point of $M(\bfx)$. For a general signomial program, items (a) and (b) are the hardest to check. Our examples provide some clues. 

The standard convergence results for the quadratic penalty method are covered in the references 
\cite{lange04,NocedalWrightBook,RuszczynskiBook}.  To summarize the principal finding, suppose that the objective function $f(\bfx)$ and the constraint functions $r_i(\bfx)$ and $s_i(\bfx)$ are continuous and that $f(\bfx)$ is coercive on $\mathbb{R}_{>0}^n$. If 
$\bfx_\lambda$ minimizes the penalized objective function 
\begin{eqnarray*}
f_{\lambda}(\bfx) & = &  f(\bfx) + \lambda \sum_i r_i(\bfx)^2 + \lambda \sum_j s_j(\bfx)_+^2 ,
\end{eqnarray*}
and $\bfx_\infty$ is a cluster point of $\bfx_\lambda$ as $\lambda$ tends to $\infty$, then $\bfx_\infty$ minimizes 
$f(\bfx)$ subject to the constraints.  In this regard observe that the coerciveness assumption on $f(\bfx)$ implies that the solution set $\{\bfx_\lambda\}_{\lambda}$ is bounded and possesses at least one cluster point. Of course, if the solution set consists of a single point, then $\bfx_{\lambda}$ tends to that point.

\section{Discussion}

The current paper presents novel algorithms for both geometric and signomial programming.  Although our examples are low dimensional, the previous experience of Sha et al.\ \cite{Sha} offers convincing evidence that the MM algorithm works well for high-dimensional quadratic programming with nonnegativity constraints.  The ideas pursued here -- the MM principle, separation of variables, quasi-Newton acceleration, and penalized optimization -- are surprisingly potent in large-scale optimization.  The MM algorithm deals with the objective function directly and reduces multivariate minimization to a sequence of one-dimensional minimizations. The MM updates are simple to code and enjoy the crucial descent property. Treating constrained signomial programming by the penalty method extends the MM algorithm even further. Quadratic programming with linear equality and inequality constraints is the most important special case of constrained signomial programming. Our new MM algorithm for constrained quadratic programming deserves consideration in high-dimensional problems. Even though MM algorithms can be notoriously slow to converge, quasi-Newton acceleration can dramatically improve matters.  Acceleration involves no matrix inversion, only matrix times vector multiplication.  Finally, it is worth keeping in mind that parameter separated algorithms are ideal candidates for parallel processing. 

Because geometric programs are ultimately convex, it is relatively easy to pose and check sufficient conditions for global convergence of the MM algorithm. In contrast it is far more difficult to analyze the behavior of the MM algorithm for signomial programs. Theoretical progress will probably be piecemeal and require problem-specific information. A major difficulty is understanding the asymptotic nature
of the objective function as parameters approach 0 or $\infty$. Even in the absence of theoretical guarantees, the descent property of the MM algorithm makes it an attractive solution technique and a diagnostic tool for finding counterexamples. Some of our test problems 
expose the behavior of the MM algorithm in non-standard situations. We welcome the help of the optimization community in unraveling the mysteries of the MM algorithm in signomial programming.




\end{document}